\documentclass[1p]{article}

\usepackage{graphicx}
\usepackage{amsmath}
\usepackage{amssymb}
\usepackage{graphicx}
\usepackage{fancybox}
\usepackage{ascmac}
\usepackage{color}
\usepackage{amsthm}
\usepackage{amscd}
\usepackage{ulem}
\usepackage[margin=20mm]{geometry}

\newtheorem{thm}{Theorem}[section]
\newtheorem{dfn}[thm]{Definition}
\newtheorem{pro}[thm]{Proposition}
\newtheorem{lemm}[thm]{Lemma}
\newtheorem{cor}[thm]{Corollary}
\newtheorem{rem}[thm]{Remark}
\newtheorem{prf}{Proof}

\begin{document}

\title{A non-tame and non-co-tame automorphism of the polynomial ring in three variables}
\author{Shoya Yasuda\footnote{Affiliation: Department of Mathematical Sciences, Tokyo Metropolitan University, 1-1 Minami-Osawa, Hachioji,
Tokyo 192-0397, Japan.}}
\date{}

\maketitle

\begin{abstract}
An automorphism $F$ of the polynomial ring in $n$ variables over a field of characteristic zero is said to be {\it co-tame} if the subgroup of the automorphism group of the polynomial ring generated by $F$ and affine automorphisms contains the tame subgroup. 
There exist many examples of such an $F$, and several sufficient conditions for co-tameness are already known. In 2015, Edo-Lewis gave the first example of a non-co-tame automorphism, which is a tame automorphism of the polynomial ring in three variables. 
In this paper, we give the first example of a non-co-tame automorphism which is not tame.
We construct such an example when $n=3$ as the exponential automorphism of a locally nilpotent derivation of rank three.
\end{abstract}

\section{Introduction}

Let $k$ be a field of characteristic zero, $k[{\bf x}]:=k[x_1,\ldots,x_n]$ the polynomial ring in $n$ variables over $k$, and $\mathrm{Aut}_kk[{\bf x}]$ the automorphism group of the $k$-algebra $k[{\bf x}]$.
We write each $\phi \in \mathrm{Aut}_kk[{\bf x}]$ as $\phi=(\phi(x_1),\ldots,\phi(x_n))$.
Given $\phi, \psi \in \mathrm{Aut}_kk[{\bf x}]$, the composition is defined by $\phi \circ \psi =(\phi(\psi(x_1)),\ldots,\phi(\psi(x_n)))$.
We say that $\phi \in \mathrm{Aut}_kk[{\bf x}]$ is {\it affine} if $\phi = (x_1,\ldots,x_n)A+\mathbf{b}$ for some $A \in \mathit{GL}_n(k)$ and $\mathbf{b} \in k^n$, and {\it triangular} if 
$$\phi = (a_1 x_1 + f_1,\ldots, a_i x_i +f_i,\ldots, a_nx_n+f_n)$$
for some $a_i \in k^*$ and $f_i \in k[x_1,\ldots,x_{i-1}]$.
We denote by $\mathrm{Aff}_n(k)$ (resp.\ $\mathrm{BA}_n(k)$) the set of affine (resp.\ triangular) automorphisms of $k[{\bf x}]$.
Then, $\mathrm{Aff}_n(k)$ and $\mathrm{BA}_n(k)$ are subgroups of $\mathrm{Aut}_kk[{\bf x}]$.

We say that $\phi \in \mathrm{Aut}_kk[{\bf x}]$ is {\it tame} if $\phi$ belongs to the {\it tame subgroup} $\mathrm{TA}_n(k):=\langle \mathrm{Aff}_n(k), \mathrm{BA}_n(k) \rangle$.
Jung \cite{Jung} and van der Kulk \cite{van1} showed that $\mathrm{Aut}_kk[x_1, x_2]=\mathrm{TA}_2(k)$.
In fact, $\mathrm{Aut}_kk[x_1, x_2]$ is the amalgamated free product of $\mathrm{Aff}_2(k)$ and $\mathrm{BA}_2(k)$ over $\mathrm{Aff}_2(k) \cap \mathrm{BA}_2(k)$ (cf.\ \cite[Part I, Theorem 3.3.]{Nagata}).
In 1972, Nagata \cite{Nagata} conjectured that $\psi \in \mathrm{Aut}_kk[x_1, x_2, x_3]$ defined by
\begin{align*}
\psi(x_1)=x_1+2(x_1x_3-x_2^2)x_2+(x_1x_3-x_2^2)^2x_3, \quad \psi(x_2)=x_2+(x_1x_3-x_2^2)x_3
\end{align*}
and $\psi(x_3)=x_3$ is not tame.
In 2004, this famous conjecture was solved in the affirmative by Shestakov-Umirbaev \cite{SU1}, \cite{SU2}.
When $n \geq 4$, it is not known whether $\mathrm{Aut}_kk[{\bf x}]=\mathrm{TA}_n(k)$.

In this paper, we study {\it co-tameness} of automorphisms defined as follows.
\begin{dfn}[Edo \cite{Edo1}] \rm
$\phi \in \mathrm{Aut}_kk[{\bf x}]$ is said to be {\it co-tame} if $\langle \phi,\mathrm{Aff}_n(k) \rangle \supset \mathrm{TA}_n(k)$.
\end{dfn}
No element of $\mathrm{Aut}_kk[x_1, x_2]$ is co-tame because of the amalgamated free product structure of $\mathrm{TA}_2(k)$ mentioned above (see \cite{EL}).
When $n \geq 3$, it is difficult to decide co-tameness of elements of $\mathrm{Aut}_kk[{\bf x}]$ in general.
The first example of a co-tame automorphism was gave by Derksen.
He showed that the triangular automorphism $(x_1,\ldots,x_{n-1},x_n+x_1^2)$ is co-tame if $n \geq 3$ (cf.\ \cite[Theorem 5.2.1]{van}).
More generally, Bodnarchuck \cite{Bod1} showed that every non-affine and one-parabolic automorphism is co-tame.
Here, $\phi \in \mathrm{Aut}_kk[{\bf x}]$ is said to be {\it one-parabolic} if there exist $a \in k^*$ and $f_1,\ldots,f_n \in k[x_1,\ldots,x_{n-1}]$ such that $\phi=(f_1,\ldots,f_{n-1},ax_n+f_n)$.
In 2013, Edo \cite{Edo1} showed that a certain class of non-tame automorphisms, including Nagata's automorphism, are co-tame.
In 2019, Edo-Lewis \cite {EL2} also gave a sufficient condition for co-tameness of automorphisms.

On the other hand, in 2015, Edo-Lewis [3] gave the first example of a non-co-tame automorphism when $n=3$. 
The example of Edo-Lewis is tame, and it is previously not known whether there exists a non-co-tame automorphism which is not tame. 

The purpose of this paper is to give the first example of a non-cotame automorphism which is not tame.
We construct such an example when $n=3$ as the exponential automorphism of a locally nilpotent derivation of rank three.
Freudenburg \cite{Freudenburg1} gave a rank three locally nilpotent derivation of $k[x_1,x_2,x_3]$.
More generally, Daigle gave a rank $n$ locally nilpotent derivation of $k[\mathbf{x}]$ when $n \geq 4$, and the corresponding exponential automorphism is one-parabolic, i.e., co-tame (cf.\ \cite[Section 3]{Freudenburg2}).
On the other hand, every exponential automorphism of a rank three locally nilpotent derivation of $k[x_1,x_2,x_3]$ is not one-parabolic because a non-zero locally nilpotent derivation of $k[x_1,x_2]$ has rank one (cf.\ \cite{Ren}).
It is remarkable that our result is the first case where co-tameness of the exponential automorphism of a rank three locally nilpotent derivation of $k[x_1,x_2,x_3]$ is decided.

\section{Main result}

A $k$-linear map $D:k[\mathbf{x}] \to k[\mathbf{x}]$ is called a {\it $k$-derivation} on $k[\mathbf{x}]$ if $D$ satisfies $D(fg)=fD(g)+D(f)g$ for all $f, g \in k[\mathbf{x}]$.
If $D$ is a $k$-derivation on $k[\mathbf{x}]$, then we can write
$$D=D(x_1)\frac{\partial}{\partial x_1} + \cdots + D(x_n) \frac{\partial}{\partial x_n}.$$
We say that a $k$-derivation $D$ on $k[\mathbf{x}]$ is {\it locally nilpotent} if, for every $f \in k[\mathbf{x}]$, there exists a positive integer $l$ such that $D^l(f)=0$.
For example, $D$ is locally nilpotent if $D(x_i)\in k[x_1,\ldots ,x_{i-1}]$ for each $i$.
Such a $k$-derivation $D$ is said to be {\it triangular}. 
We denote by $\mathrm{Der}_kk[{\bf x}]$ (resp.\ $\mathrm{LND}_kk[{\bf x}]$) the set of $k$-derivations (resp.\ locally nilpotent $k$-derivations) on $k[\mathbf{x}]$.
The {\it rank} of $D \in \mathrm{Der}_kk[{\bf x}]$ is defined to be the minimal number $r$ for which there exists $\tau \in \mathrm{Aut}_kk[{\bf x}]$ such that $D(\tau(x_i))=0$ for $i=1, \ldots, n-r$ (cf.\ \cite{Freudenburg2}).
Given $F=(f_1, \ldots , f_{n-1}) \in {k[{\bf x}]}^{n-1}$, we define the {\it Jacobian derivation} $\Delta_F \in \mathrm{Der}_kk[{\bf x}]$ by
\begin{equation*}
\Delta_F(f) := \det \frac{\partial(f_1, \ldots , f_{n-1}, f)}{\partial(x_1, \ldots , x_n)} \ \mathrm{for} \ \mathrm{all} \ f \in k[{\bf x}].
\end{equation*}
Given $D \in \mathrm{LND}_kk[{\bf x}]$, we define an {\it exponential automorphism} $\exp D \in \mathrm{Aut}_kk[{\bf x}]$ by
\begin{align*}
(\exp D)(f):=\sum_{i=0}^\infty \frac{D^i(f)}{i!} \quad \mathrm{for} \ \mathrm{all} \ f \in k[{\bf x}].
\end{align*}

In what follows, we assume that $n=3$. 
We set $f:=x_1x_3-x_2^2$, $r:=x_2f+x_1^2$ and $g:=x_3f^2+2x_1x_2f+x_1^3$, and define
$$\Delta:=\Delta_{(f,g)} \ \mathrm{and} \ \phi:=\exp \Delta.$$
It is known that $\Delta$ is a rank three locally nilpotent derivation satisfying $\Delta(r)=-fg$ and $\Delta(x_1)=-2rf$, and that $\phi$ is not tame  (cf.\ \cite{Kuroda}).
We mention that $\Delta$ is a slight modification of an example of Freudenburg \cite{Freudenburg1}.

The following is the main theorem of this section.
\begin{thm} \label{t2}
In the notation above, $\phi$ is not co-tame.
\end{thm}
Since $-fg=\Delta(r)=\Delta(x_2)f+2x_1\Delta(x_1)=\Delta(x_2)f+2x_1(-2rf)$, we have $\Delta(x_2)=4x_1r-g$.
We can also check that $\Delta(x_3)=6x_2r+2f^2$.
For each $u \in k^*$, we define $\phi_u:=\exp u \Delta$ and $\beta_u:=(u^3x_1,u^2x_2,ux_3)$, and set $B:=\{\beta_u \mid u \in k^*\}$.
Note that
\begin{align}\label{a2}
\beta_u(f)=u^4f, \ \beta_u(r)=u^6r \ \mathrm{and} \ \beta_u(g)=u^9g.
\end{align}
Using (\ref{a2}), we can check that $(\beta_u \circ \Delta)(x_i)=u^7(\Delta \circ \beta_u)(x_i)$ for $i=1,2,3$.
Hence, we get $\beta_u \circ \Delta \circ \beta_u^{-1}=u^7\Delta$.
This implies that
\begin{align}\label{a3}
\beta_{v} \circ \phi_{u} \circ \beta_{v}^{-1}=\phi_{uv^7} \ \mathrm{for \ each} \ u, v \in k^*.
\end{align}
Thanks to (\ref{a3}), we can write each $\theta \in \langle \phi, \mathrm{Aff}_3(k) \rangle \setminus \mathrm{Aff}_3(k)$ as
\begin{align}\label{a4}
\theta= \alpha_0 \circ \phi_{u_1} \circ \alpha_1 \circ \phi_{u_2} \circ \cdots \circ \phi_{u_{s-1}} \circ \alpha_{s-1} \circ \phi_{u_s} \circ \alpha_s,
\end{align}
where $s \geq 1$, $u_1,\ldots,u_s \in k^*$, $\alpha_0, \alpha_s \in \mathrm{Aff}_3(k)$ and $\alpha_1,\ldots,\alpha_{s-1} \in \mathrm{Aff}_3(k) \setminus B$.

The key to proving Theorem \ref{t2} is the following theorem.
\begin{thm} \label{t1}
For $s \geq 1$, $u_1,\ldots,u_s \in k^*$, $\alpha_0 \in \mathrm{Aff}_3(k)$ and $\alpha _1,\ldots ,\alpha _{s-1}\in \mathrm{Aff}_3(k) \setminus B$, we set $\theta = \alpha_0 \circ \phi_{u_1} \circ \alpha_1 \circ \phi_{u_2}\circ \cdots \circ \phi_{u_{s-1}} \circ \alpha_{s-1} \circ \phi_{u_s}$.
Then, we have $\deg \theta (x_i)\ge 9$ for $i=1,2,3$. 
\end{thm}
Here, ``$\deg $" denotes the total degree.
Theorem \ref{t1} implies that, for any $\theta \in \langle \phi, \mathrm{Aff}_3(k) \rangle \setminus \mathrm{Aff}_3(k)$, there exists $\alpha \in \mathrm{Aff}_3(k)$ such that $\deg(\theta \circ \alpha)(x_i) \geq 9$ for $i=1,2,3$.
Assuming Theorem \ref{t1}, we can prove Theorem \ref{t2} as follows: 
\begin{prf} [Proof of Theorem \ref{t2}] \rm
Suppose that $\phi$ is co-tame.
Then, $\langle \phi, \mathrm{Aff}_3(k) \rangle \setminus \mathrm{Aff}_3(k)$ contains the tame automorphism $\theta:=(x_1,x_2,x_3+x_1^2)$.
By the remark after Theorem \ref{t1}, there exists $\alpha \in \mathrm{Aff}_3(k)$ such that $\deg(\theta \circ \alpha)(x_1) \geq 9$, which is absurd.
\qed
\end{prf}
Theorem \ref{t1} also implies the following corollaries.
\begin{cor} \label{c2}
We have $C:=\{\alpha \in \mathrm{Aff}_3(k) \mid \alpha \circ \phi =\phi \circ \alpha \}=\{ \beta_u \in B \mid u \in k^*,  u^7=1\}$.
Hence, $C$ is a cyclic group of order seven if $k$ contains a primitive seventh root of unity, and $C=\{ \mathrm{id}_{k[\mathbf{x}]}\} $ otherwise.
\end{cor}
\begin{prf} \rm
Set $C'=\{ \beta_u \in B \mid u \in k^*, u^7=1\}$.
By (\ref{a3}) with $u=1$, we have $C \cap B=C'$.
Hence, it suffices to check that $C \subset B$.
Suppose for contradiction that there exists $\rho \in C \setminus B$.
Then, $(\phi \circ \rho \circ \phi^{-1})(x_1)=\rho(x_1)$ is a linear polynomial.
On the other hand, by Theorem \ref{t1}, we have $\deg(\phi \circ \rho \circ \phi^{-1})(x_1) \geq 9$.
This is a contradiction.
\qed
\end{prf}

By (\ref{a3}), we have $\langle \phi ,B \rangle \subset B\circ \{ \phi_u\mid u\in k\}$.
Hence, we get $\langle \phi ,B\rangle \cap \mathrm{Aff}_3(k) =B$ since $\phi_u \notin \mathrm{Aff}_3(k)$ if $u \in k^*$ (see Section \ref{sub1}, (\ref{a6})).

\begin{cor} \label{c3}
The group $\langle \phi, \mathrm{Aff}_3(k) \rangle$ is the amalgamated free product of $\mathrm{Aff}_3(k)$ and $\langle \phi, B \rangle$ over their intersection $B$.
\end{cor}
\begin{prf} \rm
Clearly, $\mathrm{Aff}_3(k)$ and $\langle \phi ,B \rangle$ generate $\langle \phi , \mathrm{Aff}_3(k)  \rangle$.
Take $\alpha_1, \ldots, \alpha_{s} \in \mathrm{Aff}_3(k) \setminus B$ and $\rho_1,\ldots,\rho_s \in \langle \phi, B \rangle \setminus B$, and set $\theta:=\alpha_1\circ \rho_1 \circ \cdots \circ \alpha_{s} \circ \rho_s$, where $s \geq 1$.
It suffices to check that $\theta \neq \mathrm{id}_{k[\mathbf{x}]}$.
For each $1 \leq i \leq s$, there exist $u_i, v_i \in k^*$ such that $\rho_i=\beta_{v_i} \circ \phi_{u_i}$, since $\langle \phi ,B \rangle \subset B\circ \{ \phi_u\mid u\in k\}$.
Then, we can write $\theta=\alpha'_1 \circ \phi_{u_1} \circ \cdots \circ \alpha'_{s} \circ \phi_{u_s}$, where $\alpha'_i := \alpha_i \circ \beta_{v_i}$.
Since $\alpha'_i \in \mathrm{Aff}_3(k) \setminus B$ for each $i$, we know by Theorem \ref{t1} that $\theta \notin \mathrm{Aff}_3(k)$.
\qed
\end{prf}

\section{Lifts of automorphisms}\label{sub1}

Now, let $n=3$.
Let $k[{\bf x},t_r,t_f,t_g]$ be the polynomial ring in six variables over $k$, and $\pi : k[{\bf x},t_r,t_f,t_g] \to k[{\bf x}]$ the substitution map defined by $t_r \mapsto r$, $t_f \mapsto f$ and $t_g \mapsto g$.
We define a triangular $k$-derivation $\Delta'$ on $k[{\bf x},t_r,t_f,t_g]$ by
\begin{align}\label{a5}
\Delta':=-2t_rt_f\frac{\partial}{\partial x_1}+(4x_1t_r-t_g)\frac{\partial}{\partial x_2}+(6x_2t_r+2t_f^2)\frac{\partial}{\partial x_3}-t_ft_g\frac{\partial}{\partial t_r}.
\end{align}
For each $u \in k^*$, we define $\phi'_u \in \mathrm{Aut}_kk[{\bf x},t_r,t_f,t_g]$ by $\phi'_u:=\exp u\Delta'$.
Then, we have $\pi \circ \phi'_u = \phi_u \circ \pi$ since $\pi \circ (u\Delta ')=(u\Delta) \circ \pi $.
By computation, we can check that
\begin{align}\label{a6}
\phi'_u(x_1)&=x_1-2ut_rt_f+u^2t_f^2t_g, \nonumber \\
\phi'_u(x_2)&=x_2+4ux_1t_r-ut_g-4u^2t_r^2t_f-2u^2x_1t_ft_g+4u^3t_rt_f^2t_g-u^4t_f^3t_g^2, \\
\phi'_u(x_3)&=x_3+6ux_2t_r+2ut_f^2-3u^2x_2t_ft_g+12u^2x_1t_r^2-3u^2t_rt_g+2u^3t_ft_g^2 \nonumber \\
&-8u^3t_r^3t_f-12u^3x_1t_rt_ft_g+3u^4x_1t_f^2t_g^2+12u^4t_r^2t_f^2t_g-6u^5t_rt_f^3t_g^2+u^6t_f^4t_g^3, \nonumber
\end{align}
$\phi'_u(t_r)=t_r-ut_ft_g$, $\phi'_u(t_f)=t_f$ and $\phi'_u(t_g)=t_g$.

For affine automorphisms, we introduce the following notion.
\begin{dfn} \rm
For $\alpha \in \mathrm{Aff}_3(k)$, we call an endomorphism $\alpha'$ of the $k$-algebra $k[\mathbf{x},t_r,t_f,t_g]$ a $\mathit{lift}$ of $\alpha$ if $\alpha'(x_i)=\alpha(x_i)$ for $i=1,2,3$ and $\pi \circ \alpha'=\alpha \circ \pi$.
\end{dfn}

Next, for $p = {\displaystyle \sum_{i_1,\ldots,i_6 \geq 0} u_{i_1,\ldots,i_6}x_1^{i_1}x_2^{i_2}x_3^{i_3}t_r^{i_4}t_f^{i_5}t_g^{i_6}} \in k[{\bf x},t_r,t_f,t_g] \setminus \{0\}$ with $u_{i_1,\ldots,i_6} \in k$, we define
$$\mathrm{supp}(p):=\{ (i_1,\ldots,i_6) \in \mathbb{N}^6 \mid u_{i_1,\ldots,i_6} \neq 0 \}.$$
For $\mathbf{w}=(w_1,\ldots,w_6) \in \mathbb{N}^6$, we define
$$\deg_{\mathbf{w}}(p):=\max\{i_1w_1+\dots+i_6w_6 \mid (i_1,\ldots,i_6) \in \mathrm{supp}(p)\}.$$
We denote by $\mathrm{lt}(p)$ the leading term of $p$ for the {\it sixth cyclic lexicographic order}, i.e., the ordering defined by $x_1^{i_1}x_2^{i_2}x_3^{i_3}t_r^{i_4}t_f^{i_5}t_g^{i_6} < x_1^{j_1}x_2^{j_2}x_3^{j_3}t_r^{j_4}t_f^{j_5}t_g^{j_6}$ if $i_6 < j_6$, or $i_6 = j_6$ and $i_m < j_m$ for the first $m$ with $i_m \neq j_m$.

Now, we set $\mathbf{w}_1:=(1,2,3,1,0,1)$ and $\mathbf{w}_2:=(0,0,0,0,1,0)$.
For each $\gamma \geq 0$ and $\delta \geq 1$, we define
\begin{align*}
\mathcal{P}_{\gamma, \delta}:=\{p \in k[\mathbf{x},t_r,t_f,t_g] \setminus \{0\} \mid \deg_{\mathbf{w}_1}(p) \leq \gamma , \deg_{\mathbf{w}_2}(p) \leq \delta, \mathrm{lt}(p) \in k^*t_f^{\delta} t_g^{\gamma}\}.
\end{align*}
From (\ref{a6}), we see that $\mathrm{lt}(\phi'_u(x_i))=(-1)^{i+1}u^{2i}t_f^{i+1}t_g^{i}$ and $\phi'_u(x_i) \in \mathcal{P}_{i,i+1}$ for $i=1,2,3$.
We remark that, if $p$ is an element of $\mathcal{P}_{\gamma, \delta}$, then one of the following holds for each $\mathbf{i}=(i_1,\ldots,i_6) \in \mathrm{supp}(p-\mathrm{lt}(p))$:

(a) $\mathbf{i}=(0,0,0,0,i_5,\gamma)$ and $0 \leq i_5 < \delta$.

(b) $0 \leq i_6 < \gamma$, $i_1+2i_2+3i_3+i_4+i_6 \leq \gamma$ and $i_5 \leq \delta$.
\begin{lemm} \label{l1}
For each $\gamma \geq 0$, $\delta \geq 1$ and $p\in \mathcal{P}_{\gamma ,\delta }$, we have $\deg \pi (p)=5\gamma +2\delta $. 
\end{lemm}
\begin{prf} \rm
Note that $d(\mathbf{i}):=\deg(\pi(x_1^{i_1}x_2^{i_2}x_3^{i_3}t_r^{i_4}t_f^{i_5}t_g^{i_6}))=\deg(x_1^{i_1}x_2^{i_2}x_3^{i_3}r^{i_4}f^{i_5}g^{i_6})=i_1+i_2+i_3+3i_4+2i_5+5i_6$ for each $\mathbf{i}=(i_1,\ldots,i_6) \in \mathrm{supp}(p)$.
If $\mathbf{i}$ is as in (a), then $d(\mathbf{i})=5\gamma+2i_5<5\gamma +2\delta$.
If $\mathbf{i}$ is as in (b), then $d(\mathbf{i}) \leq 3(i_1+2i_2+3i_3+i_4+i_6)+2i_5+2i_6 \leq 3\gamma+2\delta+2i_6<5\gamma +2\delta$.
Hence, $\deg(\pi(p))$ is equal to $\deg(\pi(\mathrm{lt}(p)))=\deg(\pi(t_f^{\delta} t_g^{\gamma}))=5\gamma +2\delta$.
\qed
\end{prf}
We prove the following proposition in the next section.
\begin{pro}\label{p1}
For any $\gamma, \delta \geq 1$, $u \in k^*$ and $\alpha \in \mathrm{Aff}_3(k) \setminus B$, there exist $\gamma' \ge \gamma$, $\delta' \geq \delta$ and a lift $\alpha'$ of $\alpha$ such that $(\phi'_u \circ \alpha') (\mathcal{P}_{\gamma, \delta}) \subset \mathcal{P}_{\gamma', \delta'}$.
\end{pro}
From Proposition \ref{p1}, we can derive Theorem \ref{t1} as follows: 
\begin{prf} [Proof of Theorem \ref{t1}] \rm
Without loss of generality, we may assume that $\alpha_0=\mathrm{id}_{k[\mathbf{x}]}$.
Fix $i=1,2,3$.
First, we show the following statement by induction on $s$:\\
There exist $\gamma \geq i$, $\delta \geq i+1$ and a lift $\alpha'_j$ of $\alpha_j$ for $j=1,\ldots,s-1$ such that
\begin{equation} \label{ghghd}
(\phi'_{u_1} \circ \alpha'_1 \circ \cdots \circ \alpha'_{s-1} \circ \phi'_{u_s})(x_i) \in \mathcal{P}_{\gamma, \delta}.
\end{equation}
Since $\phi'_{u_1}(x_i) \in \mathcal{P}_{i, i+1}$, the assertion holds when $s=1$.
Assume that $s \geq 2$.
By induction assumption, there exist $\gamma' \geq i$, $\delta' \geq i+1$ and a lift $\alpha'_j$ of $\alpha_j$ for $j=2,\ldots,s-1$ such that
$$(\phi'_{u_2} \circ \alpha'_2 \circ \cdots \circ \alpha'_{s-1} \circ \phi'_{u_s})(x_i) \in \mathcal{P}_{\gamma', \delta'}.$$
By Proposition \ref{p1}, we can find $\gamma \geq \gamma'$, $\delta \geq \delta'$ and a lift $\alpha'_1$ of $\alpha_1$ such that $(\phi_{u_1} \circ \alpha'_1)(\mathcal{P}_{\gamma', \delta'}) \subset \mathcal{P}_{\gamma, \delta}$.
Then, we have $\gamma \geq i$, $\delta \geq i+1$ and (\ref{ghghd}).

Now, since $\pi \circ \phi'_{u_i}=\phi_{u_i} \circ \pi$ and $\pi \circ \alpha'_{i}=\alpha_{i} \circ \pi$ for each $i$, we know that
\begin{align*}
\theta(x_i)&=(\phi_{u_1} \circ \alpha_1 \circ \cdots \circ \alpha_{s-1} \circ \phi_{u_s} \circ \pi)(x_i)\\
&=(\pi \circ \phi'_{u_1} \circ \alpha'_1 \circ \cdots \circ \alpha'_{s-1} \circ \phi'_{u_s})(x_i) \in \pi(\mathcal{P}_{\gamma, \delta})
\end{align*}
by (\ref{ghghd}).
This implies that $\deg (\theta (x_i)) = 5\gamma +2\delta \geq 5i +2(i+1) \geq 9$ thanks to Lemma \ref{l1}.
\qed
\end{prf}

\section{Proof of Proposition \ref{p1}}

The rest of the paper is devoted to the proof of Proposition \ref{p1}.
Fix $\gamma, \delta \geq 1$, $u\in k^*$ and $\alpha \in \mathrm{Aff}_3(k) \setminus B$.
Write 
\begin{align*}
\alpha=(x_1,x_2,x_3)\left(
\begin{array}{rrr}
a_{1,1} & a_{1,2} & a_{1,3} \\
a_{2,1} & a_{2,2} & a_{2,3} \\
a_{3,1} & a_{3,2} & a_{3,3} \\
\end{array}
\right)
+(d_1,d_2,d_3)
\end{align*}
and set $\gamma_j:=\max\{i \mid a_{i,j} \neq 0\}$ for $j=1,2,3$.
We call $(\gamma_1, \gamma_2, \gamma_3)$ the $\mathit{type}$ of $\alpha$.
Since $(a_{i,j}) \in \mathit{GL}_3(k)$, we have $\gamma_i=3$ for some $i$, and $\gamma_j \geq 2$ for some $j \neq i$.
For $i=1,2,3$, we define 
$$X_i:=\phi'_u(\alpha(x_i))=a_{1,i}\phi'_u(x_1)+a_{2,i}\phi'_u(x_2)+a_{3,i}\phi'_u(x_3)+d_i.$$
Since $\phi'_u(x_j) \in \mathcal{P}_{j, j+1}$ for $j=1,2,3$ as mentioned, we see that $\mathrm{lt}(X_i)=a_{\gamma_i,i}\mathrm{lt}(\phi'_u(x_{\gamma_i}))$ and $X_i \in \mathcal{P}_{\gamma_i, \gamma_i+1}$.
We note that $X_i=\phi'_u(\alpha'(x_i))$ for any lift $\alpha'$ of $\alpha$.

For each lift $\alpha'$ of $\alpha$, we consider the following condition (C).\\
(C) There exist integers $m_i \geq 0$ and $n_i \geq 1$ for $i=4,5,6$ such that\\
(c1) $\phi'_u(\alpha'(t_r)) \in \mathcal{P}_{m_4, n_4}$, $\phi'_u(\alpha'(t_f)) \in \mathcal{P}_{m_5, n_5}$ and $\phi'_u(\alpha'(t_g)) \in \mathcal{P}_{m_6, n_6}$;\\
(c2) $m_6 \geq \gamma_1$, $2m_6 \geq \gamma_2$, $3m_6 \geq \gamma_3$ and $m_6 \geq m_4$;\\
(c3) $n_6-1 \geq \gamma_1+1$, $2(n_6-1) \geq \gamma_2+1$, $3(n_6-1) \geq \gamma_3+1$  and $n_6-1 \geq n_4$.

\begin{lemm}\label{lele}
Let $\alpha'$ be a lift of $\alpha$ which satisfies $\mathrm{(C)}$, and set $\gamma':=m_5\delta+m_6\gamma$ and $\delta':=n_5\delta+n_6\gamma$.
Then, we have $\gamma' \geq \gamma$, $\delta' \geq \delta$ and $(\phi'_u \circ \alpha')(\mathcal{P}_{\gamma, \delta}) \subset \mathcal{P}_{\gamma', \delta'}$.
\end{lemm}

\begin{prf} \rm
Since $m_6 \geq \gamma_1\geq 1$ by (c2), and $n_5 \geq 1$, we have $\gamma' \geq \gamma$ and $\delta' \geq \delta$.
Set $\psi:=\phi'_u \circ \alpha'$.
Take any $p \in \mathcal{P}_{\gamma, \delta}$.
For $\mathbf{i}=(i_1,\ldots,i_6) \in \mathrm{supp}(p)$ and $\mathbf{w} \in \{\mathbf{w}_1, \mathbf{w}_2 \}$, we set
$$d_{\mathbf{i}, \mathbf{w}} := \deg_{\mathbf{w}}(\psi(x_1^{i_1}x_2^{i_2}x_3^{i_3}t_r^{i_4}t_f^{i_5}t_g^{i_6}))=\deg_{\mathbf{w}}(X_1^{i_1}X_2^{i_2}X_3^{i_3}\psi(t_r)^{i_4}\psi(t_f)^{i_5}\psi(t_g)^{i_6}).$$
By the definition of $\mathcal{P}_{\gamma, \delta}$, we have $i_5 \leq \delta$, $i_1+2i_2+3i_3+i_4+i_6 \leq \gamma$ and $i_6 \leq \gamma$. 
Noting $X_i \in \mathcal{P}_{\gamma_i, \gamma_i+1}$ for $i=1,2,3$ and (c1), we get
\begin{align*}
d_{\mathbf{i}, \mathbf{w}_1}&\leq \gamma_1i_1+\gamma_2 i_2+\gamma_3 i_3+m_4 i_4+m_5 i_5+m_6 i_6 \\
&\leq m_6 i_1+2m_6 i_2 +3m_6 i_3 +m_6 i_4+m_5 i_5+m_6 i_6 & (\mathrm{\because \ (c2)})\\
&= m_5 i_5+m_6(i_1+2i_2+3i_3+i_4+i_6) \leq m_5 \delta + m_6 \gamma,\\
d_{\mathbf{i}, \mathbf{w}_2}&\leq (\gamma_1+1) i_1+(\gamma_2+1) i_2+(\gamma_3+1) i_3+n_4 i_4+n_5 i_5+n_6 i_6 \\
&\leq (n_6-1)i_1+2(n_6-1)i_2+3(n_6-1)i_3+(n_6-1)i_4+n_5i_5+n_6 i_6 & (\mathrm{\because \ (c3)})\\
&= n_5 i_5+(n_6-1) (i_1+2i_2+3i_3+i_4+i_6) +i_6 \\
&\leq n_5 \delta +(n_6-1)\gamma+i_6 \leq n_5 \delta +n_6 \gamma.\\
\end{align*}
This proves that $\deg_{\mathbf{w}_1}(\psi(p)) \leq m_5 \delta + m_6 \gamma$ and $\deg_{\mathbf{w}_2}(\psi(p)) \leq n_5 \delta +n_6 \gamma$.
Next, we consider $\mathrm{lt}(\psi(p))$.
Noting $X_i \in \mathcal{P}_{\gamma_i, \gamma_i+1}$ for $i=1,2,3$ and (c1), we have $\mathrm{lt}(\psi(x_1^{i_1}x_2^{i_2}x_3^{i_3}t_r^{i_4}t_f^{i_5}t_g^{i_6}))=\mathrm{lt}(X_1^{i_1}X_2^{i_2}X_3^{i_3}\psi(t_r)^{i_4}\psi(t_f)^{i_5}\psi(t_g)^{i_6}) \in k^*m_{\mathbf{i}}$, where
$$m_{\mathbf{i}} := t_f^{(\gamma_1+1) i_1+(\gamma_2+1) i_2+(\gamma_3+1) i_3+n_4 i_4+n_5 i_5+n_6 i_6}t_g^{\gamma_1i_1+\gamma_2 i_2+\gamma_3 i_3+m_4 i_4+m_5 i_5+m_6 i_6}.$$
If $\mathbf{i}$ is as in (a), then $m_{\mathbf{i}}=t_f^{n_5i_5+n_6\gamma}t_g^{m_5 i_5+m_6\gamma}$ is less than $t_f^{n_5 \delta +n_6 \gamma}t_g^{m_5 \delta +m_6 \gamma}$, since $i_5 < \delta$.
If $\mathbf{i}$ is as in (b), then $(\gamma_1+1) i_1+(\gamma_2+1) i_2+(\gamma_3+1) i_3+n_4 i_4+n_5 i_5+n_6 i_6 \leq n_5 \delta +(n_6-1)\gamma+i_6 < n_5 \delta + n_6 \gamma $ since $i_6 < \gamma$.
Hence, $m_{\mathbf{i}}$ is less than $t_f^{n_5 \delta +n_6 \gamma}t_g^{m_5 \delta +m_6 \gamma}$ for the sixth cyclic lexicographic order.
Thus, we conclude that $\mathrm{lt}(\psi(p)) \in k^*t_f^{n_5 \delta +n_6 \gamma}t_g^{m_5 \delta +m_6 \gamma}$.
Therefore, $\psi(p)$ belongs to $\mathcal{P}_{m_5 \delta + m_6 \gamma, n_5 \delta +n_6 \gamma}$.
\qed
\end{prf}

We define an endomorphism $\hat{\alpha}$ of the $k$-algebra $k[{\bf x},t_r,t_f,t_g]$ by $\hat{\alpha}(x_i):=\alpha(x_i)$ for each $i=1,2,3$, and
\begin{align*}
\hat{\alpha}(t_f):&=\hat{\alpha}(f)+(a_1c_3-2a_2c_2+a_3c_1)(t_f-f),\\
\hat{\alpha}(t_r):&=\hat{\alpha}(x_2t_f+x_1^2),\\
\hat{\alpha}(t_g):&=\hat{\alpha}(x_3t_f^2+2x_1x_2t_f+x_1^3).
\end{align*}
It is clear that $(\pi \circ \hat{\alpha})(x_i)=(\alpha \circ \pi)(x_i)$ for each $i=1,2,3$.
Since $\pi(t_f-f)=0$, we have $(\pi \circ \hat{\alpha})(t_f)=\pi(\hat{\alpha}(f))=\alpha(f)=(\alpha \circ \pi)(t_f)$.
Thus, we get
\begin{align*}
(\pi \circ \hat{\alpha})(t_r)=(\pi \circ \hat{\alpha})(x_2t_f+x_1^2)=(\alpha \circ \pi)(x_2t_f+x_1^2)=\alpha(x_2f+x_1^2)=(\alpha \circ \pi)(t_r).
\end{align*}
Similarly, we can check that $(\pi \circ \hat{\alpha})(t_g)=(\alpha \circ \pi)(t_g)$.
Therefore, $\hat{\alpha}$ is a lift of $\alpha$.
We will often use this $\hat{\alpha}$ in the proof of Proposition \ref{p1}.
So, we put $\psi:=\phi'_u \circ \hat{\alpha}$, and $R:=\psi(t_r)$, $F:=\psi(t_f)$ and $G:=\psi(t_g)$.

\begin{rem} \label{rem000} \rm
Let $\alpha'$ be an endomorphism of the $k$-algebra $k[\mathbf{x},t_r,t_f,t_g]$ with $\alpha'(x_i)=\alpha(x_i)$ for $i=1,2,3$.
If $\alpha'(t_{\lambda})-\hat{\alpha}(t_{\lambda}) \in \ker \pi$ for all $\lambda \in \{r,f,g\}$, then $\alpha'$ is also a lift of $\alpha$.
\end{rem}

Now, we prove Proposition \ref{p1}.
Thanks to Lemma \ref{lele}, it suffices to construct a lift $\hat{\alpha}$ of $\alpha$ which satisfies (C).
We consider the following five cases separately.

\begin{lemm}\label{l2}
If $(\gamma_1,\gamma_3)=(3,3)$ and $\gamma_2 \leq 2$, or $(\gamma_1,\gamma_3) \neq (3,3)$ and $\gamma_2=3$, then $\hat{\alpha}$ satisfies $\mathrm{(C)}$.
\end{lemm}

\begin{prf} \rm
For simplicity, we write $a_j:=a_{1,j}$, $b_j:=a_{2,j}$ and $c_j:=a_{3,j}$ for $j=1,2,3$.
Then, we have
\begin{align*}
\hat{\alpha}(t_f)&=(c_1c_3-c_2^2)x_3^2+(b_1c_3-2b_2c_2+b_3c_1)x_2x_3+(a_1c_3-2a_2c_2+a_3c_1+b_1b_3-b_2^2)x_2^2\\
& \quad +(a_1b_3-2a_2b_2+a_3b_1)x_1x_2+Px_3+(a_1a_3-a_2^2)x_1^2+Qx_2+Sx_1\\
& \quad +(a_1c_3-2a_2c_2+a_3c_1)t_f+d_1d_3-d_2^2,
\end{align*}
where $P:=c_1d_3-2c_2d_2+c_3d_1$, $Q:=b_1d_3-2b_2d_2+b_3d_1$ and $S:=a_1d_3-2a_2d_2+a_3d_1$.
If $(\gamma_1,\gamma_3)=(3,3)$ and $\gamma_2 \leq 2$, then $c_1c_3-c_2^2=c_1c_3 \neq 0$.
If $(\gamma_1,\gamma_3) \neq (3,3)$ and $\gamma_2=3$, then $c_1c_3-c_2^2=-c_2^2 \neq 0$.
Thus, the monomial $x_3^2$ appears in $\alpha '(t_f)$ with nonzero coefficient. 
Since $\phi'_u(x_3) \in \mathcal{P}_{3,4}$, it follows that $F \in \mathcal{P}_{6,8}$.

Now, $G=\psi(x_3t_f^2+2x_1x_2t_f+x_1^3)=X_3F^2+2X_1X_2F+X_1^3$.
Since $F \in \mathcal{P}_{6,8}$, we have $X_3F^2 \in \mathcal{P}_{12+\gamma_3,17+\gamma_3}$, $2X_1X_2F \in \mathcal{P}_{6+\gamma_1+\gamma_2,10+\gamma_1+\gamma_2}$ and $X_1^3 \in \mathcal{P}_{3\gamma_1,3+3\gamma_1}$.
Since $1 \leq \gamma_i \leq 3$ for each $i=1,2,3$, we know that $G \in \mathcal{P}_{12+\gamma_3,17+\gamma_3}$.
Similarly, $R=X_2F+X_1^2$ belongs to $\mathcal{P}_{6+\gamma_2,9+\gamma_2}$.
In this case, (c1), (c2) and (c3) hold for $(m_4,m_5,m_6)=(6+\gamma_2,6,12+\gamma_3)$ and $(n_4,n_5,n_6)=(9+\gamma_2,8,17+\gamma_3)$.
Therefore, $\hat{\alpha}$ satisfies (C).
\qed
\end{prf}

We prove the following lemmas similarly.

\begin{lemm}\label{l3}
If $(\gamma_1,\gamma_3) \in \{(3,2), (2,3)\}$ and $\gamma_2 \leq 2$, then $\hat{\alpha}$ satisfies $\mathrm{(C)}$.
\end{lemm}

\begin{prf} \rm
If $(\gamma_1,\gamma_3)=(3,2)$ and $\gamma_2 \leq 2$, then $c_1c_3-c_2^2=0$ and $b_1c_3-2b_2c_2+b_3c_1=b_3c_1\neq0$.
If $(\gamma_1,\gamma_3)=(2,3)$ and $\gamma_2 \leq 2$, then $c_1c_3-c_2^2=0$ and $b_1c_3-2b_2c_2+b_3c_1=b_1c_3\neq0$.
Thus, the monomial $x_2x_3$ appears in $\alpha '(t_f)$ with nonzero coefficient.
This implies that $F \in \mathcal{P}_{5,7}$.
Hence, we have $X_3F^2 \in \mathcal{P}_{10+\gamma_3,15+\gamma_3}$, $2X_1X_2F \in \mathcal{P}_{5+\gamma_1+\gamma_2,9+\gamma_1+\gamma_2}$ and $X_1^3 \in \mathcal{P}_{3\gamma_1,3+3\gamma_1}$.
Since $1 \leq \gamma_i \leq 3$ for each $i=1,2,3$, it follows that $G \in \mathcal{P}_{10+\gamma_3,15+\gamma_3}$.
Similarly, $R$ belongs to $\mathcal{P}_{5+\gamma_2,8+\gamma_2}$.
Then, we can easily check that $\hat{\alpha}$ satisfies (C).
\qed
\end{prf}

We remark that $\hat{\alpha}(t_f) \notin k$, since $\pi(\hat{\alpha}(t_f))=\alpha(\pi(t_f))=\alpha(f) \notin k$.

\begin{lemm}\label{l4}
If $(\gamma_1,\gamma_2,\gamma_3)=(3,2,1)$, then $\hat{\alpha}$ satisfies $\mathrm{(C)}$.
\end{lemm}
\begin{prf} \rm
Since $(\gamma_1,\gamma_2,\gamma_3)=(3,2,1)$, we have $X_1 \in \mathcal{P}_{3,4}$, $X_2 \in \mathcal{P}_{2,3}$ and $X_3 \in \mathcal{P}_{1,2}$, and $c_1 \neq 0$, $b_2 \neq 0$, $a_3 \neq 0$ and $c_2=c_3=b_3=0$.
Hence, we get
\begin{align*}
\hat{\alpha}(t_f)&=(a_3c_1-b_2^2)x_2^2+(-2a_2b_2+a_3b_1)x_1x_2+c_1d_3x_3+(a_1a_3-a_2^2)x_1^2\\
& \quad +(b_1d_3-2b_2d_2)x_2+(a_1d_3-2a_2d_2+a_3d_1)x_1+a_3c_1t_f+d_1d_3-d_2^2.
\end{align*}
We remark that, if $F \in \mathcal{P}_{\gamma,\delta}$ for some $\gamma \geq 0$ and $\delta \geq 1$, then $X_3F^2 \in \mathcal{P}_{1+2\gamma,2+2\delta}$, $2X_1X_2F \in \mathcal{P}_{5+\gamma,7+\delta}$ and $X_1^3 \in \mathcal{P}_{9,12}$, and $X_2F \in \mathcal{P}_{2+\gamma,3+\delta}$ and $X_1^2 \in \mathcal{P}_{6,8}$.

Suppose that  $a_3c_1 - b_2^2 \neq 0$.
Then, since the monomial $x_2^2$ appears in $\alpha '(t_f)$ with nonzero coefficient, we have $F \in \mathcal{P}_{4,6}$.
By the remark, we get $G \in \mathcal{P}_{9,14}$ and $R \in \mathcal{P}_{6,9}$.
In this case, $\hat{\alpha}$ satisfies (C).

Suppose that  $a_3c_1 - b_2^2 = 0$.
Then, since $\hat{\alpha}(t_f) \notin k$ as remarked, there exists $$(\gamma, \delta) \in \{(0,1), (1,2), (2,3), (2,4), (3,4), (3,5)\}$$ such that $F \in \mathcal{P}_{\gamma, \delta}$.
Since $\gamma \leq 3$ and $\delta \leq 5$, we get $G \in \mathcal{P}_{9,12}$ and $R \in \mathcal{P}_{6,8}$ by the remark.
In this case, $\hat{\alpha}$ satisfies (C).
\qed
\end{prf}

\begin{lemm}\label{l5}
If $(\gamma_1,\gamma_2,\gamma_3)=(1,2,3)$, then there exists a lift of $\alpha$ which  satisfies $\mathrm{(C)}$.
\end{lemm}
\begin{prf} \rm
Since $(\gamma_1,\gamma_2,\gamma_3)=(1,2,3)$, we have $X_1 \in \mathcal{P}_{1,2}$, $X_2 \in \mathcal{P}_{2,3}$ and $X_3 \in \mathcal{P}_{3,4}$, and $a_1 \neq 0$, $b_2 \neq 0$, $c_3 \neq 0$ and $c_1=c_2=b_1=0$.
Hence, we get
\begin{align*}
\hat{\alpha}(t_f)&=(a_1c_3- b_2^2)x_2^2+(a_1b_3 - 2a_2b_2)x_1x_2+c_3d_1x_3+(a_1a_3 - a_2^2)x_1^2\\
& \quad +(- 2b_2d_2+ b_3d_1)x_2+(a_1d_3 - 2a_2d_2+ a_3d_1)x_1+a_1c_3t_f+d_1d_3 - d_2^2.
\end{align*}
We remark that, if $F \in \mathcal{P}_{\gamma,\delta}$ for some $\gamma \geq 0$ and $\delta \geq 1$, then $X_3F^2 \in \mathcal{P}_{3+2\gamma,4+2\delta}$, $2X_1X_2F \in \mathcal{P}_{3+\gamma,5+\delta}$ and $X_1^3 \in \mathcal{P}_{3,6}$, and $X_2F \in \mathcal{P}_{2+\gamma,3+\delta}$ and $X_1^2 \in \mathcal{P}_{2,4}$.

\begin{flushleft}
(i) $a_1c_3- b_2^2 \neq 0$, $a_1b_3 - 2a_2b_2 \neq 0$, $d_1 \neq 0$, $a_1a_3 - a_2^2 \neq 0$, $d_2 \neq 0$, or $d_3 \neq 0$.
\end{flushleft}

In this case, there exists $$(\gamma, \delta) \in \{(1,2), (2,3), (2,4), (3,4), (3,5), (4,6)\}$$ such that $F \in \mathcal{P}_{\gamma, \delta}$.
Then, since $\gamma \geq 1$ and $\delta \geq 2$, we get $G \in \mathcal{P}_{3+2\gamma,4+2\delta}$ and $R \in \mathcal{P}_{2+\gamma,3+\delta}$ by the remark.
In this case, $\hat{\alpha}$ satisfies (C).

In the following, we assume that $a_1c_3-b_2^2=a_1b_3 - 2a_2b_2=d_1=a_1a_3 - a_2^2=d_2=d_3=0$.
Then, we have $\hat{\alpha}(t_f)=a_1c_3t_f$, so $F=\phi'_u(\hat{\alpha}(t_f))=a_1c_3t_f \in \mathcal{P}_{0,1}$.
By the remark, it follows that $X_3F^2, 2X_1X_2F, X_1^3 \in \mathcal{P}_{3,6}$ and $X_2F, X_1^2 \in \mathcal{P}_{2,4}$.
\begin{flushleft}
(ii) $b_2-c_3^2\neq0$.
\end{flushleft}

Since $\mathrm{lt}(X_3F^2)=u^6a_1^2c_3^3t_f^6t_g^3$, $\mathrm{lt}(2X_1X_2F)=-2u^6a_1^2b_2c_3t_f^6t_g^3$ and $\mathrm{lt}(X_1^3)=u^6a_1^3t_f^6t_g^3$, we have
\begin{align*}
\mathrm{lt}(X_3F^2)+\mathrm{lt}(2X_1X_2F)+\mathrm{lt}(X_1^3)=u^6a_1^2(a_1-2b_2c_3+c_3^3)t_f^6t_g^3.
\end{align*}
Since $\mathrm{lt}(X_2F)=-u^4a_1b_2c_3t_f^4t_g^2$ and $\mathrm{lt}(X_1^2)=u^4a_1^2t_f^4t_g^2$, we have
\begin{align*}
\mathrm{lt}(X_2F)+\mathrm{lt}(X_1^2)=u^4a_1(a_1-b_2c_3)t_f^4t_g^2.
\end{align*}
Now, since $a_1c_3-b_2^2 = 0$, we have $\displaystyle a_1 - 2b_2c_3 + c_3^3=\frac{1}{c_3}(b_2-c_3^2)^2 \neq 0$ and $\displaystyle a_1 - b_2c_3=\frac{b_2}{c_3}(b_2-c_3^2) \neq 0$.
Thus, we get $G \in \mathcal{P}_{3,6}$ and $R \in \mathcal{P}_{2,4}$.
In this case, $\hat{\alpha}$ satisfies (C).

\begin{flushleft}
(iii) $b_2-c_3^2=0$ and $a_2 \neq 0$.
\end{flushleft}

Since $a_1c_3- b_2^2 = 0$ and $b_2-c_3^2=0$, we have $a_1=c_3^3$ and $b_2=c_3^2$.
Since $a_1b_3 - 2a_2b_2=0$, it follows that $2a_2=b_3c_3$.
Hence, we can write $\alpha(x_1)=c_3^3x_1$, $\alpha(x_2)=(b_3c_3/2)x_1+c_3^2x_2$ and $\alpha(x_3)=a_3x_1+b_3x_2+c_3x_3$, and
\begin{align*}
\hat{\alpha}(t_f)&=c_3^4t_f,\\
\hat{\alpha}(t_g)&=(a_3x_1+b_3x_2+c_3x_3)(c_3^4t_f)^2+2(c_3^3x_1)((b_3c_3/2)x_1+c_3^2x_2)(c_3^4t_f)+(c_3^3x_1)^3\\
&=a_3c_3^8x_1t_f^2 + b_3c_3^8t_f(x_2t_f + x_1^2) + c_3^9(x_3t_f^2 + 2x_1x_2t_f + x_1^3),\\
\hat{\alpha}(t_r)&=((b_3c_3/2)x_1+c_3^2x_2)(c_3^4t_f)+(c_3^3x_1)^2\\
&=(b_3c_3^5/2)x_1t_f+c_3^6(x_2t_f+x_1^2).
\end{align*}
By Remark \ref{rem000}, we can define a lift $\alpha'$ of $\alpha$ by $\alpha'(t_f):=\hat{\alpha}(t_f)$,
\begin{align*}
\alpha'(t_g)&:=\hat{\alpha}(t_g) +b_3c_3^8t_f(t_r-(x_2t_f + x_1^2))+ c_3^9(t_g-(x_3t_f^2 + 2x_1x_2t_f + x_1^3))\\
&=a_3c_3^8x_1t_f^2 + b_3c_3^8t_rt_f +c_3^9t_g,\\
\alpha'(t_r)&:=\hat{\alpha}(t_r) + c_3^6(t_r-(x_2t_f + x_1^2))=(b_3c_3^5/2)x_1t_f + c_3^6t_r.
\end{align*}

Since $a_1a_3 - a_2^2 = 0$, $2a_2-b_3c_3=0$ and $a_1, c_3 \neq 0$, we have $a_2 \neq 0 \Leftrightarrow a_3 \neq 0 \Leftrightarrow b_3 \neq 0$.
Since $\phi'_u(t_r)=t_r-ut_ft_g \in \mathcal{P}_{1,1}$, we see that $(\phi'_u \circ \alpha')(t_r)=a_2c_3^4X_1t_f+c_3^6\phi'_u(t_r) \in \mathcal{P}_{1,3}$, $(\phi'_u \circ \alpha')(t_f)=c_3^4t_f \in \mathcal{P}_{0,1}$ and $(\phi'_u \circ \alpha')(t_g)=a_3c_3^8X_1t_f^2 + b_3c_3^8\phi'_u(t_r)t_f +c_3^9t_g \in \mathcal{P}_{1,4}$.
In this case, $\alpha'$ satisfies (C).

\begin{flushleft}
(iv) If $b_2-c_3^2=a_2= 0$, then $\alpha \in B$, which contradicts the choice of $\alpha$. \qed
\end{flushleft}
\end{prf}

\begin{lemm}\label{l6}
If $\gamma_1=\gamma_2=\gamma_3=3$, then  there exists a lift of $\alpha$ which satisfies $\mathrm{(C)}$.
\end{lemm}
\begin{prf} \rm
Since $\gamma_1=\gamma_2=\gamma_3=3$, we have $X_i \in \mathcal{P}_{3,4}$ and $c_i \neq 0$ for $i=1,2,3$.
We remark that, if  $F \in \mathcal{P}_{\gamma,\delta}$ for some $\gamma \geq 0$ and $\delta \geq 1$, then $X_3F^2 \in \mathcal{P}_{3+2\gamma,4+2\delta}$, $2X_1X_2F \in \mathcal{P}_{6+\gamma,8+\delta}$ and $X_1^3 \in \mathcal{P}_{9,12}$, and $X_2F \in \mathcal{P}_{3+\gamma,4+\delta}$ and $X_1^2 \in \mathcal{P}_{6,8}$.

\begin{flushleft}
(i) $c_1c_3-c_2^2 \neq 0$, $b_1c_3-2b_2c_2+b_3c_1 \neq 0$, $a_1c_3-2a_2c_2+a_3c_1+b_1b_3-b_2^2 \neq 0$, or $a_1b_3-2a_2b_2+a_3b_1 \neq 0$.
\end{flushleft}

In this case, there exists  $$(\gamma, \delta) \in \{(3,5), (4,6), (5,7), (6,8) \}$$ such that $F \in \mathcal{P}_{\gamma, \delta}$.
Then, since $\gamma \geq 3$ and $\delta \geq 5$, we get $G \in \mathcal{P}_{3+2\gamma,4+2\delta}$ and $R \in \mathcal{P}_{3+\gamma,4+\delta}$ by the remark.
In this case, $\hat{\alpha}$ satisfies (C). 

In the following, we assume that $c_1c_3-c_2^2=b_1c_3-2b_2c_2+b_3c_1=a_1c_3-2a_2c_2+a_3c_1+b_1b_3-b_2^2=a_1b_3-2a_2b_2+a_3b_1=0$.

\begin{flushleft}
(ii) $P=c_1d_3-2c_2d_2+c_3d_1 = 0$.
\end{flushleft}

In this case, there exists $$(\gamma, \delta) \in \{(0,1), (1,2), (2,3), (2,4) \}$$ such that $F \in \mathcal{P}_{\gamma, \delta}$, since $\hat{\alpha}(t_f) \notin k$.
Then, since $\gamma \leq 2$ and $\delta \leq 4$, we get $G \in \mathcal{P}_{9,12}$ and $R \in \mathcal{P}_{6,8}$ by the remark.
In this case, $\hat{\alpha}$ satisfies (C).

In the following, we assume that $P \neq 0$.
Then, we have $F \in \mathcal{P}_{3,4}$, and so $X_3F^2,2X_1X_2F,X_1^3 \in \mathcal{P}_{9,12}$ and $X_2F, X_1^2 \in \mathcal{P}_{6,8}$ by the remark. 

\begin{flushleft}
(iii) $Pc_2+c_1^2 \neq 0$.
\end{flushleft}

Since $\mathrm{lt}(X_3F^2)=u^{18}P^2c_3t_f^{12}t_g^9$, $\mathrm{lt}(2X_1X_2F)=2u^{18}Pc_1c_2t_f^{12}t_g^9$ and $\mathrm{lt}(X_1^3)=u^{18}c_1^3t_f^{12}t_g^9$, we have
\begin{align*}
\mathrm{lt}(X_3F^2)+\mathrm{lt}(2X_1X_2F)+\mathrm{lt}(X_1^3)=u^{18}(P^2c_3+2Pc_1c_2+c_1^3)t_f^{12}t_g^9.
\end{align*}
Since $\mathrm{lt}(X_2F)=u^{12}Pc_2t_f^8t_g^6$ and $\mathrm{lt}(X_1^2)=u^{12}c_1^2t_f^8t_g^6$, we have
\begin{align*}
\mathrm{lt}(X_2F)+\mathrm{lt}(X_1^2)=u^{12}(Pc_2+c_1^2)t_f^8t_g^6.
\end{align*}
Now, since $c_1c_3-c_2^2 = 0$, we have $\displaystyle P^2c_3+2Pc_1c_2+c_1^3=\frac{1}{c_1}(Pc_2+c_1^2)^2 \neq 0$.
Thus, we get $G \in \mathcal{P}_{9,12}$ and $R \in \mathcal{P}_{6,8}$.
In this case, $\hat{\alpha}$ satisfies (C). 

In the following, we assume that $Pc_2+c_1^2 = 0$.
Since we already assumed that $P\neq0$, $c_1c_3-c_2^2=0$, $b_1c_3 - 2b_2c_2 + b_3c_1=0$ and $a_1c_3-2a_2c_2+a_3c_1+b_1b_3-b_2^2=0$, we get
\begin{gather}
(c_1,c_2,c_3)=\left(c_1,-\frac{1}{P}c_1^2,\frac{1}{P^2}c_1^3\right),\label{rete} \\
b_3=-\frac{1}{P^2}b_1c_1^2-\frac{2}{P}b_2c_1 \quad \mathrm{and} \quad a_3=\frac{1}{P^2c_1}(P^2b_2^2 -2Pa_2c_1^2 + 2Pb_1b_2c_1 - a_1c_1^3 + b_1^2c_1^2).\label{rete2}
\end{gather}
We also remark that
\begin{equation}\label{8989}
\begin{aligned}
P^2d_3 + 2Pc_1d_2 + c_1^2d_1&=\frac{P^2}{c_1}\left(c_1d_3 + 2\frac{1}{P}c_1^2d_2 + \frac{1}{P^2}c_1^3d_1\right)\\
&=\frac{P^2}{c_1}(c_1d_3 - 2c_2d_2 +c_3d_1)\\
&=\frac{P^3}{c_1}\neq0.
\end{aligned}
\end{equation}
By (\ref{rete}) and the first equality of (\ref{rete2}), we have
\begin{equation}\label{eee1}
\begin{aligned}
X_1&=\phi'_u(a_1x_1+b_1x_2+c_1x_3+d_1),\\
X_2&=\phi'_u\left(a_2x_1+b_2x_2-\frac{1}{P}c_1^2x_3+d_2\right),\\
X_3&=\phi'_u\left(a_3x_1+\left(-\frac{1}{P^2}b_1c_1^2-\frac{2}{P}b_2c_1\right)x_2+\frac{1}{P^2}c_1^3x_3+d_3\right).
\end{aligned}
\end{equation}
By (\ref{rete}) and the second equality of (\ref{rete2}), we have $a_1c_3-2a_2c_2+a_3c_1=(1/P^2)(Pb_2+b_1c_1)^2$.
Thus, we get
\begin{equation}\label{eee2}
\begin{aligned}
F=\phi'_u\left(Px_3+(a_1a_3-a_2^2)x_1^2+Qx_2+Sx_1+\frac{1}{P^2}(Pb_2+b_1c_1)^2t_f+d_1d_3-d_2^2\right).
\end{aligned}
\end{equation}

\begin{flushleft}
(iv) $a_1a_3-a_2^2 \neq 0$. 
\end{flushleft}

Using (\ref{a6}), (\ref{eee1}) and (\ref{eee2}), we can describe $G=X_3F^2+2X_1X_2F+X_1^3$ explicitly, and check that
$$\deg_{\mathbf{w}_1}(G)=7, \ \deg_{\mathbf{w}_2}(G)=12 \ \mathrm{and} \ \mathrm{lt}(G)=\frac{1}{P^2}u^{14}c_1^3(a_1a_3-a_2^2)^2t_f^{12}t_g^7$$ 
by computer.
Thus, we get $G \in \mathcal{P}_{7,12}$.
Similarly, we can describe $R=X_2F+X_1^2$ explicitly, and check that
$$\deg_{\mathbf{w}_1}(R)=5, \ \deg_{\mathbf{w}_2}(R)=8 \ \mathrm{and} \ \mathrm{lt}(R)=-\frac{1}{P}u^{10}c_1^2(a_1a_3 - a_2^2)t_f^8t_g^5$$
by computer.
Thus, we get $R \in \mathcal{P}_{5,8}$.
In this case, $\hat{\alpha}$ satisfies (C). 

In the following, we assume that $a_1a_3-a_2^2 = 0$.
Then, (\ref{eee2}) implies
\begin{equation}\label{eee4}
F=\phi'_u\left(Px_3+Qx_2+Sx_1+\frac{1}{P^2}(Pb_2+b_1c_1)^2t_f+d_1d_3-d_2^2\right).
\end{equation}

\begin{flushleft}
(v) $P^2b_2 + 2Pb_1c_1 - Qc_1^2 \neq 0$. 
\end{flushleft}

Using (\ref{a6}), (\ref{eee4}) and (\ref{eee1}) with $a_3$ replaced by the second part of (\ref{rete2}), we can describe $G$ explicitly, and check that
$$\deg_{\mathbf{w}_1}(G)=7, \ \deg_{\mathbf{w}_2}(G)=10 \ \mathrm{and} \ \mathrm{lt}(G)=\frac{1}{P^2c_1}u^{14}(P^2b_2 + 2Pb_1c_1 - Qc_1^2)^2t_f^{10}t_g^7$$
by computer.
Thus, we get $G \in \mathcal{P}_{7,10}$.
Similarly, we can describe $R$ explicitly, and check that
$$\deg_{\mathbf{w}_1}(R)=5, \ \deg_{\mathbf{w}_2}(R)=7 \ \mathrm{and} \ \mathrm{lt}(R)=-\frac{1}{P}u^{10}(P^2b_2 + 2Pb_1c_1 - Qc_1^2)t_f^7t_g^5$$
by computer.
Thus, we get $R \in \mathcal{P}_{5,7}$.
In this case, $\hat{\alpha}$ satisfies (C). 

In the following, we assume that $P^2b_2 + 2Pb_1c_1 - Qc_1^2 = 0$.
Set $T:=Pb_1-Qc_1$.
Then, we get $P^2b_2+Pb_1c_1+c_1T=0$.
From this equality and the first equality of (\ref{rete2}), it follows that
\begin{equation}\label{gav}
(b_1,b_2,b_3)=\left(b_1,-\frac{T}{P^2}c_1-\frac{1}{P}b_1c_1,2\frac{T}{P^3}c_1^2+\frac{1}{P^2}b_1c_1^2\right).
\end{equation}
We remark that $T \neq 0$.
Actually, if $T=0$, then $b_1=(Q/P)c_1$, so we have $(b_1,b_2,b_3)=(Q/P)(c_1,-(1/P)c_1^2,(1/P^2)c_1^3) \in k \cdot (c_1,c_2,c_3)$ by (\ref{rete}).
This contradicts $(a_{i,j}) \in \mathit{GL}_3(k)$.
We show that
\begin{equation}\label{hgu}
(a_1,a_2,a_3)=\left(\frac{1}{4c_1}b_1^2,-\frac{T}{2P^2}b_1-\frac{1}{4P}b_1^2,\frac{T^2}{P^4}c_1+\frac{T}{P^3}b_1c_1+\frac{1}{4P^2}b_1^2c_1\right).
\end{equation}
Substituting the second component of (\ref{gav}) for $b_2$ in the second equality of (\ref{rete2}), we get
\begin{equation}\label{gavv}
a_3=\frac{T^2}{P^4}c_1-2\frac{1}{P}a_2c_1-\frac{1}{P^2}a_1c_1^2.
\end{equation}
By $(\ref{gav})$ and $(\ref{gavv})$, we have $0=a_1b_3-2a_2b_2+a_3b_1=(Tc_1/P^4)(2Pa_1c_1+2P^2a_2+Tb_1)$.
This gives 
\begin{equation}\label{gavvv}
a_2=-\frac{1}{P}a_1c_1-\frac{T}{2P^2}b_1.
\end{equation}
By $(\ref{gavv})$ and $(\ref{gavvv})$, we have $0=a_1a_3-a_2^2=(T^2/4P^4)(4a_1c_1-b_1^2)$.
This gives $a_1=(1/4c_1)b_1^2$.
From this equality and (\ref{gavv}) and (\ref{gavvv}), we get $(\ref{hgu})$.

Hereafter, we consider a lift $\alpha'$ of $\alpha$ defined by $\alpha'(t_f):=\hat{\alpha}(t_f)$,
\begin{equation*}
\begin{aligned}
\alpha'(t_g)&:=\hat{\alpha}(t_g)+\frac{T^2}{P^2}c_1x_3(t_f-(x_1x_3-x_2^2)),\\
\alpha'(t_r)&:=\hat{\alpha}(t_r)-\frac{T^2}{P^2}(t_f-(x_1x_3-x_2^2)).
\end{aligned}
\end{equation*}
In this case, $(\phi'_u \circ \alpha')(t_f)$ is the same as $F$.
We set $G':=(\phi_u' \circ \alpha')(t_g)$ and $R':=(\phi_u' \circ \alpha')(t_r)$.

By (\ref{eee1}) with (\ref{gav}) and $(\ref{hgu})$, we have
\begin{equation}\label{eee5}
\begin{aligned}
X_1&=\phi'_u\left(\frac{1}{4c_1}b_1^2x_1+b_1x_2+c_1x_3+d_1\right),\\
X_2&=\phi'_u\left(\left(-\frac{T}{2P^2}b_1-\frac{1}{4P}b_1^2\right)x_1+\left(-\frac{T}{P^2}c_1-\frac{1}{P}b_1c_1\right)x_2-\frac{1}{P}c_1^2x_3+d_2\right),\\
X_3&=\phi'_u\left(\left(\frac{T^2}{P^4}c_1+\frac{T}{P^3}b_1c_1+\frac{1}{4P^2}b_1^2c_1\right)x_1+\left(2\frac{T}{P^3}c_1^2+\frac{1}{P^2}b_1c_1^2\right)x_2+\frac{1}{P^2}c_1^3x_3+d_3\right).
\end{aligned}
\end{equation}
Since $Q=(1/c_1)(Pb_1-T)$ by the definition of $T$, we have
\begin{equation}\label{eee6}
F=\phi'_u\left(Px_3+\frac{1}{c_1}(Pb_1-T)x_2+Sx_1+\frac{T^2}{P^4}c_1^2t_f+d_1d_3-d_2^2\right)
\end{equation}
by (\ref{eee4}), (\ref{gav}).

Using (\ref{a6}), (\ref{eee5}) and (\ref{eee6}), we can describe $G'$ explicitly, and check that
$$\deg_{\mathbf{w}_1}(G')=6, \ \deg_{\mathbf{w}_2}(G')=8 \ \mathrm{and} \ \mathrm{lt}(G')=u^{12}(P^2d_3 + 2Pc_1d_2 + c_1^2d_1)t_f^8t_g^6$$
by computer, where we note that $P^2d_3 + 2Pc_1d_2 + c_1^2d_1\neq0$ by (\ref{8989}).
Thus, $G'$ belongs to $\mathcal{P}_{6,8}$.
We shall show that $R'$ belongs to one of $\mathcal{P}_{4,6}$, $\mathcal{P}_{3,5}$ and $\mathcal{P}_{3,4}$.
Then, we know that $\alpha'$ satisfies (C).

\begin{flushleft}
(vi) $P^2b_1^2 - 4PSc_1^2 - 2 PTb_1 + 4T^2 \neq 0$.
\end{flushleft}

Using (\ref{a6}), (\ref{eee5}) and (\ref{eee6}), we can describe $R'$ explicitly, and check that
$$\deg_{\mathbf{w}_1}(R')=4, \ \deg_{\mathbf{w}_2}(R')=6 \ \mathrm{and} \ \mathrm{lt}(R')=\frac{1}{4P^2}u^8(P^2b_1^2 - 4PSc_1^2 - 2 PTb_1 + 4T^2)t_f^6t_g^4$$
by computer. 
Thus, we get $R' \in \mathcal{P}_{4,6}$.

In the following, we assume that $P^2b_1^2 - 4PSc_1^2 - 2 PTb_1 + 4T^2 = 0$, that is,
\begin{equation}\label{hshs}
S=\frac{1}{4Pc_1^2}(P^2b_1^2  - 2 PTb_1 + 4T^2).
\end{equation}
Then, by (\ref{eee6}), we have
\begin{equation}\label{eee8}
\begin{aligned}
F=\phi'_u\left(Px_3+\frac{1}{c_1}(Pb_1-T)x_2+\frac{1}{4Pc_1^2}(P^2b_1^2  - 2 PTb_1 + 4T^2)x_1+\frac{T^2}{P^4}c_1^2t_f+d_1d_3-d_2^2\right).
\end{aligned}
\end{equation}

\begin{flushleft}
(vii) $P^2T - c_1^5 \neq 0$. 
\end{flushleft}

Using (\ref{a6}), (\ref{eee5}) and (\ref{eee8}), we can describe $R'$ explicitly, and check that
$$\deg_{\mathbf{w}_1}(R')=3, \ \deg_{\mathbf{w}_2}(R')=5 \ \mathrm{and} \ \mathrm{lt}(R')=\frac{1}{P^5c_1}u^6T^2(P^2T - c_1^5)t_f^5t_g^3$$ 
by computer.
Thus, we get $R' \in \mathcal{P}_{3,5}$.

\begin{flushleft}
(viii) $P^2T - c_1^5 = 0$, that is, $T=c_1^5/P^2$.
\end{flushleft}

First, we show that
\begin{equation}\label{hshshs}
(d_1,d_2,d_3)=\left(\frac{P^3}{c_1^3},-\frac{3P^2}{2c_1^2},\frac{3P}{c_1}\right).
\end{equation}
Note that
\begin{align}
c_1d_3-2c_2d_2+c_3d_1&=P,\label{111} \\
b_1d_3-2b_2d_2+b_3d_1&=Q=\frac{1}{c_1}Pb_1-\frac{1}{P^2}c_1^4,\label{222} \\
a_1d_3-2a_2d_2+a_3d_1&=S=\frac{1}{4c_1^2}Pb_1^2-\frac{1}{2P^2}b_1c_1^3+\frac{1}{P^5}c_1^8,\label{333}
\end{align}
where the second equalities of (\ref{222}) and (\ref{333}) are due to $Q=(1/c_1)(Pb_1-T)$ and (\ref{hshs}) with $T=c_1^5/P^2$, respectively.
Computing both sides of $(-P^5/c_1^5)(b_1(\ref{111})-c_1(\ref{222}))$, we get
\begin{equation}\label{444}
2Pc_1^2d_2+2c_1^3d_1=-P^3,
\end{equation}
where we use (\ref{rete}) and (\ref{gav}) for the left-hand side.
Similarly, from $(-P^8/c_1^4)(b_1(\ref{222})-4c_1(\ref{333}))$, we get
\begin{equation}\label{555}
2P^4b_1c_1^2d_2+(4c_1^8+2P^3b_1c_1^3)d_1=4P^3c_1^5-P^6b_1,
\end{equation}
where we use (\ref{gav}) and (\ref{hgu}) for the left-hand side.
Then, $(\ref{555})-P^3b_1(\ref{444})$ gives $4c_1^8d_1=4P^3c_1^5$.
This implies that $d_1=P^3/c_1^3$.
Therefore, we obtain (\ref{hshshs}) by (\ref{444}), (\ref{111}) and (\ref{rete}).

By (\ref{eee5}) with (\ref{hshshs}) and $T=c_1^5/P^2$, we have
\begin{equation}\label{eee9}
\begin{aligned}
X_1&=\phi'_u\left(\frac{1}{4c_1}b_1^2x_1+b_1x_2+c_1x_3+\frac{P^3}{c_1^3}\right),\\
X_2&=\phi'_u\left(\left(-\frac{1}{2P^4}b_1c_1^5-\frac{1}{4P}b_1^2\right)x_1+\left(-\frac{1}{P^4}c_1^6-\frac{1}{P}b_1c_1\right)x_2-\frac{1}{P}c_1^2x_3-\frac{3P^2}{2c_1^2}\right),\\
X_3&=\phi'_u\left(\left(\frac{1}{P^8}c_1^{11}+\frac{1}{P^5}b_1c_1^6+\frac{1}{4P^2}b_1^2c_1\right)x_1+\left(2\frac{1}{P^5}c_1^7+\frac{1}{P^2}b_1c_1^2\right)x_2+\frac{1}{P^2}c_1^3x_3+\frac{3P}{c_1}\right).
\end{aligned}
\end{equation}
By (\ref{eee8}) with (\ref{hshshs}) and $T=c_1^5/P^2$, we have
\begin{equation}\label{eee10}
\begin{aligned}
F=\phi'_u\left(Px_3+\frac{1}{c_1}\left(Pb_1-\frac{c_1^5}{P^2}\right)x_2+\frac{1}{4c_1^2}\left(Pb_1^2  -\frac{2}{P^2}b_1c_1^5 + \frac{4}{P^5}c_1^{10}\right)x_1+\frac{1}{P^8}c_1^{12}t_f+\frac{3P^4}{4c_1^4}\right).
\end{aligned}
\end{equation}

Using (\ref{a6}), (\ref{eee9}) and (\ref{eee10}), we can describe $R'$ explicitly, and check that
$$\deg_{\mathbf{w}_1}(R')=3, \ \deg_{\mathbf{w}_2}(R')=4 \ \mathrm{and} \ \mathrm{lt}(R')=-\frac{P^3}{4c_1^2}u^6 t_f^4t_g^3$$
by computer.
Thus, we get $R' \in \mathcal{P}_{3,4}$.
\qed
\end{prf}

This completes the proof of Proposition \ref{p1}.

\section*{Acknowledgment}
The author thanks Professor Shigeru Kuroda for useful discussions.
This work was supported by JSPS KAKENHI Grant Number JP19J20334.

\end{document}